\documentclass[12pt]{amsart}
\usepackage{epic,eepic,latexsym, amssymb, amscd, amsfonts, xypic}
\input xy
\xyoption{all}
\providecommand{\mathscr}

 \newlength{\baseunit}               % the basic unit length
 \newcount{\numlines}                % depth of picture (in number of lines)
\newcommand{\point}{\vspace{3mm}\par\refstepcounter{subsection}\noindent{\bf \thesubsection.} }
\newcommand{\tpoint}[1]{\vspace{3mm}\par\refstepcounter{subsection}\noindent{\bf \thesubsection.} 
  {\em #1. ---} }
\newcommand{\epoint}[1]{\vspace{3mm}\par\refstepcounter{subsection}\noindent{\bf \thesubsection.} 
  {\em #1.} }
\newcommand{\bpoint}[1]{\vspace{3mm}\par\refstepcounter{subsection}\noindent{\bf \thesubsection.} 
  {\bf #1.} }
\newcommand{\bpf}{\noindent {\em Proof.  }}
\newcommand{\epf}{\qed \vspace{+10pt}}

\newcommand{\Z}{\mathbb{Z}}
\newcommand{\A}{\mathbb{A}}
\newcommand{\Q}{\mathbb{Q}}

\newcommand{\C}{\mathbb{C}}
\newcommand{\proj}{\mathbb P}
\newcommand{\oh}{{\mathcal{O}}}

\newcommand{\Mbar}{\overline{M}}
\newcommand{\ch}{{\mathcal{H}}}

\newcommand{\de}{\delta}

\newcommand{\codim}{\operatorname{codim}}
\newcommand{\dimsup}{\operatorname{d}} % dimension of support of a sheaf.
\newcommand{\et}{{\mbox{\scriptsize{\'et}}}}

\newcommand{\acn}{\operatorname{acn}}
\newcommand{\easn}{\operatorname{easn}}
\newcommand{\asn}{\operatorname{asn}}
\newcommand{\clspt}[1]{\overline{\{#1\}}} % macro for the closure of a point
\newcommand{\cd}{\operatorname{cd}}
\newcommand{\psv}{\operatorname{psv}}
\newcommand{\cF}{{\mathcal{F}}}
\newcommand{\tsF}{\mathsf{F}}

\newcommand{\Fsh}{\cF}

\newcommand{\PP}{\proj}

\newcommand{\Pid}{\mathcal{P}}

\newcommand{\lremind}[1]{{}}
\newcommand{\bremind}[1]{{}}

\newcommand{\cut}[1]{}

\begin{document}
\pagestyle{plain} \title{{ \large{ 
The affine stratification number
and the moduli space of curves}}}
\author{Mike Roth}
\address{Dept.\ of Mathematics and Statistics, Queens University, Kingston,
Ontario, Canada}
\email{mikeroth@mast.queensu.ca}
\author{Ravi Vakil}
\address{Dept.\ of Mathematics, Stanford University, Stanford, CA, USA}
\email{vakil@math.stanford.edu}
\thanks{The first author is partially supported by NSERC. 
The second author is partially supported by  NSF grant DMS--0228011 and NSF
CAREER grant DMS--0238532. \newline
 }
\date{Friday, June 18, 2004.}
\begin{abstract}
  We define the {\em affine stratification number} $\asn X$ of a
  scheme $X$.  For $X$ equidimensional, it is the minimal number $k$
  such that there is a stratification of $X$ by locally closed
  affine subschemes of codimension at most $k$.  We show that the
  affine stratification number is well-behaved, and bounds many
  aspects of the topological complexity of the scheme, such as
  vanishing of cohomology groups of quasicoherent, constructible, and
  $\ell$-adic sheaves.  We explain how to bound $\asn X$ in practice.
  We give a series of conjectures (the first by E. Looijenga) bounding
  the affine stratification number of various moduli spaces of pointed curves.
  For example, the philosophy of \cite[Theorem~$\star$]{gv} yields:
  the moduli space of genus $g$, $n$-pointed complex curves of compact
  type (resp.\ with ``rational tails'') should have 
  the homotopy type of a
  finite complex of dimension at most $5g-6+2n$ (resp.\ $4g-5+2n$).
  This investigation is based on work and questions of Looijenga.

  One relevant example (Example~\ref{nosmooth}) turns out to be a proper
  integral variety with no embeddings in a smooth algebraic space.
  This one-paragraph construction appears to be simpler and more
  elementary than the earlier examples, due to Horrocks \cite{ho} and
  Nori \cite{n}.
\end{abstract}
\maketitle
\tableofcontents

{\parskip=12pt % closing bracket is just before the bibliography
\section{Introduction}
\label{intro}
The {\em affine stratification} number of a scheme $X$ bounds the
``topological complexity'' of a scheme.  For example, it bounds the
{\em cohomological dimension} $\cd X$ of $X$, which is the largest
integer $n$ such that $H^n(X, \cF) \neq 0$ for some quasicoherent
sheaf $\cF$ (Proposition~\ref{cdleqasn}).  Similarly, the cohomology of any
constructible or $\ell$-adic sheaf vanishes in degree greater than
$\asn X + \dim X$ (Proposition~\ref{fly}).
We expect that if the base field is $\C$, then $X$ has the homotopy
type of a finite complex of dimension at most $\asn X + \dim X$
(Conjecture~\ref{homtype}), but have not completed a proof.
(Unless otherwise stated, all schemes and stacks
are assumed to be separated and of finite type over an arbitrary base
field.)

A related, previously studied invariant is the {\em affine covering
  number} $\acn X$, which is one less than the minimal number of
affine open sets required to cover $X$.  The affine stratification
number is bounded by $\acn X$, is better behaved ({\em e.g.}\ is bounded by dimension, {\em cf.}\ Example~\ref{jason}), and has the same
topological consequences.  We know of no interesting consequences of
bounded $\acn$ that are not already consequences of the same bound
on $\asn$.

For equidimensional $X$, the definition is particularly simple.

\epoint{Definition} The (equidimensional) affine stratification number
of an equidimensional scheme $X$ is the minimal number $\easn X$ such
that there is a (finite) stratification of $X$ by locally closed affine
subschemes of codimension at most $\easn X$.\label{defeasn}\lremind{easn}

This is the form most likely to be of interest.  The
appropriate generalization to arbitrary schemes is only slightly more complicated.

\epoint{Definition} \label{defasn}\lremind{defasn}An {\em affine
  stratification} of a scheme $X$ is a finite decomposition $
X = \bigsqcup_{k \in \Z^{\geq 0}, i} Y_{k,i} $ into disjoint locally closed
affine subschemes $Y_{k,i}$, where for each $Y_{k,i}$,\lremind{bob}
\begin{equation}\label{bob}
\overline{Y}_{k,i}\setminus Y_{k,i} \, \subseteq \,\bigcup_{
  k'>k, \, j} Y_{k',j}.\end{equation}
The {\em length} of an affine stratification is the largest $k$ such that 
$\cup_j Y_{k,j}$ is nonempty.
The {\em affine stratification number} $\asn X$ of a scheme $X$
is the minimum of the length over all possible affine
stratifications of $X$.  

The
inclusion in \eqref{bob} refers to the underlying set.  We do not require
that each $Y_{k,i}$ be irreducible.  We also do not require any
relation between $k$ and the dimension or codimension of $Y_{k,i}$ in
$X$.  We will see however (Theorem~\ref{thm:stratathm}) that it is
always possible to assume that the stratification has a very nice
form.  

Strictly speaking, the term ``stratification'' is inappropriate, as
$\overline{Y}_{k,i}\setminus Y_{k,i}$ need not be a union of
$Y_{k',j}$: let $X$ be the co-ordinate axes in $\A^2$, $Y_{0,1}$ the
$x$-axis  minus the origin, and $Y_{1,1}$ the $y$-axis.  However,
Theorem~\ref{thm:stratathm}(a) shows that we may take \eqref{bob} to
be an actual stratification.

The affine stratification number has many good properties, including
the following (Lemma~\ref{triviallemma}, Propositions~\ref{asnleqdim},
\ref{asnleqacn}, \ref{macbeth}).
\begin{itemize}
\item $\asn X=0$ if and only if $X$ is affine.
\item $\asn X \leq \dim X$.  
\item  $\asn X \leq \acn X$.  (Equality does not always hold.)
\item $\asn (X \times Y) \leq \asn X + \asn Y.$ 
\item If $D$ is an effective Cartier divisor on $X$, then $\asn (X-D) \leq \asn X$.
\item If $Y \rightarrow X$ is an affine morphism, then $\asn Y \leq \asn X$.
\end{itemize}
Even if one is only interested in equidimensional schemes, the
more general Definition~\ref{defasn} has advantages over
Definition~\ref{defeasn}.  For example, the last property
is immediate using Definition~\ref{defasn}, but not 
obvious using Definition~\ref{defeasn}.

In Section~\ref{basic}, we establish basic properties of affine
stratifications.  In Section~\ref{basic2}, we show that affine
stratifications can be reorganized into a particularly good form.  In
particular, if $X$ is equidimensional, then $\easn X =
\asn X$ (Proposition~\ref{equidimX}), so the notation $\easn$ may be 
discarded.  In
Section~\ref{topcon}, we give topological consequences of bounded
$\asn$.

Our motivation is to bound the affine stratification number of moduli
spaces (in particular, of pointed curves) to obtain topological and
cohomological consequences.  We describe our work in progress in the
form of several conjectures in Section~\ref{ms}. For example, the
conjectures bound the homotopy type of the moduli spaces of curves (a)
of compact type (stable curves whose dual graph is a tree, or
equivalently stable curves with compact Jacobian), (b) with ``rational
tails'' (stable genus $g$ curves having a smooth component of genus
$g$), and (c) with at most $k$ rational components (a locus introduced
in \cite{gv}), see Proposition~\ref{homMgn}.

%The proper integral threefold with no smooth embeddings,
%promised in the abstract, is Example~\ref{nosmooth}.

\epoint{Acknowledgments} This note arose from our ongoing efforts to
prove a conjecture of E. Looijenga, and much of what is here derives
from questions, ideas, and work of his.  In particular, we suspect
that he is aware of most of the results given here, and that we are
following in his footsteps.  We thank him for inspiration.  We also
thank J. Starr for Example~\ref{jason}, W. Fulton for pointing out the
examples of Nori \cite{n} and Horrocks \cite{ho}, and T. Graber for
helpful conversations.  Finally, we are grateful to the organizers of
the 2003 conference {\em Algebraic structures and moduli spaces}, at
the Centre de Recherches Math\'{e}matiques (CRM), which led to this
work.

\section{Basic properties of affine stratifications} 
\label{basic}\lremind{basic}The 
most basic property is that an affine stratification always exists,
and hence $\asn X$ is defined for any scheme $X$: 
if $\cup_{i=0}^{n} U_i$ is a covering of $X$ by open affine sets, then
\lremind{exists}
\begin{equation}
\label{exists}
U_0 \sqcup (U_1 \setminus U_0) \sqcup (U_2 \setminus (U_0 \cup U_1))
\sqcup \cdots 
\end{equation}
gives an affine stratification of $X$.

The following lemma is trivial.

\tpoint{Lemma} {\em 
\begin{enumerate}
\item[(a)]  The affine stratification number depends only on the
reduced structure of $X$,\\
{\em i.e.}\ $\asn X = \asn X^{\text{red}}$.
\item[(b)]  If $X \rightarrow Y$ is an affine morphism, then
$\asn X \leq \asn Y$.  
\item[(c)]  $\asn (X \times Y) \leq \asn X + \asn Y$.
\item[(d)]  If $D$ is an effective Cartier divisor on $X$, then 
$\asn (X-D) \leq \asn X$.
\end{enumerate}} 

\label{triviallemma}\lremind{triviallemma}Part (d) requires the
following well-known fact.

\tpoint{Lemma} {\em
Any irreducible affine scheme $X$, minus an effective
Cartier divisor $D$, is affine.}
 \lremind{minuscartier}\label{minuscartier}

(Reason:  the inclusion $X \setminus D \hookrightarrow X$ is
an affine morphism, since this can be verified locally.  But $X$ is affine.)

Here is a partial converse to Lemma~\ref{minuscartier}.  
A more precise converse is given in
Proposition~\ref{affinevanishingforlcgroups}.

\tpoint{Lemma} {\em Suppose that $V$ is an irreducible affine scheme,
  and that $U\subset V$ is an open affine subset.  Then the complement
  $Z:=V\setminus U$ is a Weil divisor in $V$.
}\label{codimone}\lremind{codimone}

\bpf We first assume that $V$ (and hence $U$) is normal.  Let
$Z=\cup_i Z_i$ be the decomposition of $Z$ into irreducible components, and let
$Z'=\cup_j Z_j$ be the union of those components of
codimension one in $V$.  We set $U'=V\setminus Z'$, and let
$i:U\hookrightarrow U'$ be the natural open immersion.  Since $U'$ is
normal, and the complement of $U$ in $U'$ is of codimension at least $2$ in
$U'$, we have $i_{*}\oh_U=\oh_{U'}$.  We will use this and the fact
that both $U$ and $V$ are affine to see that $U=U'$.

Let $A=\Gamma(V,\oh_V)$ and $B=\Gamma(U,\oh_U)=\Gamma(U',\oh_{U'})$.
We have an inclusion of rings $A\hookrightarrow B$ corresponding to
the opposite inclusion of open sets.  Suppose that $U\neq U'$, and let
$x$ be any point of $U'\setminus U$.  Since $V$ is affine, $x$
corresponds to a prime ideal $\Pid_x$ of $A$.  Since $x\in U'$, no
element of $\Pid_x$ can be a unit in $\Gamma(U',\oh_{U'})$, and hence
$\Pid_x$ remains a prime ideal in $B$, which is a localization of $A$.
Therefore, since $U$ is affine, $x\in U$, contrary to assumption.

Passing to the general case, we drop the assumption that $V$ and $U$ are normal, and let $\widetilde{V}$ and $\widetilde{U}$ be their normalizations.
We have the commutative diagram
$$
\xymatrix{
\widetilde{U}\, \ar@{^{(}->}[r] \ar[d] & \widetilde{V}\ar[d] \\
U\, \ar@{^{(}->}[r]  & V  }
$$
where the vertical arrows are the normalization maps, and the
horizontal arrows are open immersions.  By the first part of the
lemma, the complement $\widetilde{Z}$ of $\widetilde{U}$ in
$\widetilde{V}$ is of codimension one in $\widetilde{V}$.  Since
$\widetilde{Z}$ maps finitely and surjectively onto $Z$,
$\dim(Z)=\dim(\widetilde{Z})$, and hence $Z$ is of codimension one in
$V$. \epf

The next corollary follows immediately.  (Note that $X$ need
not be equidimensional here.)

\tpoint{Corollary} {\em The complement of a dense affine open subset
  in any scheme is of pure codimension one.
}\label{cor:purecodimone}\lremind{cor:purecodimone}

\epoint{Examples} {\em (a)} Let $X$ be the affine cone over an
elliptic curve, embedded in $\C \PP^2$ as a cubic.  Let $Z$ be the
cone over any point of the curve of infinite order in the group law.
Then $X\setminus Z$ is affine, but $Z$ is not $\Q$-Cartier.  This
shows that the complement of an affine open set in an affine scheme
need not be the support of a Cartier divisor: we cannot hope to
improve the conclusion of Lemma~\ref{codimone} to match the hypothesis
of Lemma~\ref{minuscartier}.

{\em (b)}  Let $S$ be $\PP^2$ blown up at a point, and let $X$
be the affine cone over some projective embedding of $S$.  
Let $Z\subset X$ be the affine cone
over the exceptional divisor of the blowup.  Then $Z$ is of
codimension one in $X$, but $\cd (X\setminus Z)=1$, so in particular it
is not affine.  This shows that, conversely, the complement of a Weil divisor in
an affine scheme need not be affine:  we can not hope to
improve the hypothesis of Lemma~\ref{minuscartier} to match the
conclusion of Lemma~\ref{codimone}. 

However, there is a more precise statement giving
a necessary and sufficient condition on a closed subset
$Z$ of an affine scheme $V$ for the complement $V\setminus Z$ to be
affine.

\tpoint{Proposition} {\em Let $V$ be an affine scheme (possibly reducible)
  and $Z$ a closed subset of $V$.  Then $U:=V\setminus Z$ is affine if
  and only if $H_Z^i(\cF)=0$ for all quasicoherent sheaves $\cF$ on
  $V$ and all $i\geq 2$.}
\label{affinevanishingforlcgroups}\lremind{affinevanishingforlcgroups}

Here $H_Z^i(\cF)$ is the local cohomology group.  This also implies
the same fact for the local cohomology sheaves $\ch_Z^i(\cF)$, see
Corollary~\ref{hip}(a) below.

\bpf
Let $\cF$ be any quasicoherent sheaf on $V$.  We have the long exact excision sequence of cohomology groups\lremind{excision}
\begin{equation}
\label{excision}
\begin{array}{ccccccccc}
0 & \longrightarrow & H_Z^0(\cF) & \longrightarrow & H^0(V,\cF) & \longrightarrow & H^0(U,\cF|_U) & \\
 & \longrightarrow & H_Z^1(\cF) & \longrightarrow & H^1(V,\cF) & \longrightarrow & H^1(U,\cF|_U) & \\
 & \longrightarrow & H_Z^2(\cF) & \longrightarrow & H^2(V,\cF) & \longrightarrow & H^2(U,\cF|_U) & \longrightarrow & \cdots . \\
\end{array}
\end{equation}
Since $V$ is affine, we have $H^i(V,\cF)=0$ for all $i\geq 1$, 
so that $H^i(U,\cF|_U)=H^{i+1}_Z(\cF)$ for all $i\geq 1$.
Hence (using Serre's criterion for affineness)
 $U$ is affine if and only if $H^i_Z(\cF)=0$ for all $i\geq 2$ and all 
quasicoherent
sheaves $\cF$. \epf

\tpoint{Corollary} 
{\em Let $X$ be a scheme (possibly reducible) and
  $U$ a dense affine open subset.  Let $Z:=X\setminus U$.  
For
  any quasicoherent sheaf $\cF$ on $X$,  \label{hip}\lremind{hip}
\begin{enumerate}
\item[(a)]  $\ch^i_Z(\cF)=0$ for all
  $i\geq 2$, and
\item[(b)] $H_Z^i(\cF)=0$ for all $i>\cd Z +1$.
\end{enumerate}
}

The notation $\cd$ denotes cohomological dimension, see Section \ref{intro}.

\bpf
{\em (a)}
Since the local cohomology sheaf $\ch^i_Z(\cF)$ is the sheafification of the 
functor $V\mapsto H^i_{Z\cap V}(\cF|_V)$ \cite[Proposition\ 1.2]{lc}, 
it is sufficient
to check that the local cohomology group vanishes for sufficiently small $V$ 
around any point of $Z$.   But if $V$ is any open affine set, then $V\cap U$ 
is nonempty
(since $U$ is dense) and also affine (by separatedness). Hence
$H^i_{Z\cap V}(\cF|_V)=0$ for $i\geq 2$ by 
Proposition~\ref{affinevanishingforlcgroups} and so $\ch^i_Z(\cF)=0$ as well.

{\em (b)} The local cohomology sheaves $\ch^i_Z(\cF)$ are quasicoherent and
are supported on $Z$.  The local cohomology groups can be computed by
a spectral sequence with $E_2^{pq}$ term
$H^p(X,\ch_Z^q(\cF))=H^p(Z,\ch_Z^q(\cF))$.  Since $H^p(Z,\cdot)=0$ for
$p>\cd Z$, and $\ch^q_Z(\cF)=0$ for $q>1$ by
part (a), we have $H_Z^{i}(\cF)=0$ for
$i>\cd Z+1$. \epf

\tpoint{Corollary} {\em Let $X$ be a scheme, $U$ a dense affine open subset, and 
set $Z:=X\setminus U$.  Then $\cd X \leq \cd Z +1$. }
\label{inductiveboundingstep}\lremind{inductiveboundingstep}

\bpf
For any any quasicoherent sheaf $\cF$ on $X$, the excision sequence
\eqref{excision} and the fact that $U$ is affine gives
$H^i(X,\cF)=H_Z^{i}(\cF)$ for all $i\geq 2$, and that $H^1(X,\cF)$ is a quotient of $H^1_Z(\cF)$.  Hence, for any
$i\geq 1$, $H^i_Z(\cF)=0$ implies that $H^i(X,\cF)=0$.  Since $H^i_Z(\cF)=0$ for all $i>\cd Z+1$ by Corollary~\ref{hip}(b),
we have $\cd X\leq \cd Z+1$. \epf

\bpoint{Bounding $\asn$ by finite flat covers}
The following result is useful to bound $\asn X$ by studying
covers of $X$.

\tpoint{Proposition} {\em Suppose $\pi:Y \rightarrow X$ is a
  surjective finite flat morphism of degree not divisible by the
  characteristic of the base field, and $Y$ is affine.  Then $X$ is
  affine.}
\label{macbeth}\lremind{macbeth}

\bpf The hypothesis implies that $\pi_{*}\oh_Y$ is a coherent locally
free sheaf on $X$.  The trace map gives a splitting $\pi_{*}\oh_Y
\cong \oh_X \oplus E$ for some vector bundle $E$ on $X$.
If $\cF$ is any coherent sheaf on $X$, then the flatness of $\pi$ gives
$\pi_{*}\pi^{*}\cF=\cF \oplus (E\otimes\cF)$, and it then follows from
the Leray spectral sequence and 
the finiteness of $\pi$ that $H^i(X,\cF)$ is a direct summand of
$H^i(Y,\pi^{*}\cF)$ for all $i\geq 0$.  Since $Y$ is affine, these vanish
if $i\geq 1$, hence the cohomology groups on $X$ do as well, and therefore
$X$ is affine by Serre's criterion for affineness.
\epf

\section{Reorganizing affine stratifications}
\label{basic2}\lremind{basic2}We 
describe various ways that we can reorganize the stratification which
are more convenient for analyzing $X$.
The main results of this section are summarized
in the following theorem.

\tpoint{Theorem}{\em If $X$ is any scheme and $\asn X =m$, then there exists
an affine stratification $\{Z_0,\ldots,Z_{m}\}$ of $X$ such that for any $k\leq m$:
\begin{itemize}
\item[(a)] $\overline{Z}_k = \cup_{k'\geq k} Z_{k'},$
\item[(b)] each $Z_k$ is a dense open affine subset of $\overline{Z}_k$, and
\item[(c)] $\overline{Z}_{k}$ is of pure codimension 
one in $\overline{Z}_{k-1}$ .
\end{itemize}

\noindent
If in addition $X$ is equidimensional, then we also have
\begin{itemize}
\item[(d)] each $\overline{Z}_{k'}$ is of pure codimension $k'-k$ 
in $\overline{Z}_k$ for
any $k'\geq k$.
\end{itemize}
\noindent
Even if $X$ is not equidimensional, if we have an affine stratification $\{Y'_{k,l}\}$ of
length $M$ such that each $Y'_{k,i}$ is of pure codimension $k$ in $X$, 
then setting $Z'_k:=\cup_{i} Y'_{k,i}$ for $k=0,\ldots, M$ we have
\begin{itemize}
\item[(e)] the affine stratification $\{Z'_0,\ldots,Z'_{M}\}$ satisfies (a)--(d) above.
\end{itemize}
(We are not guaranteed that $M=m$, so this stratification
may not be optimal.)
}\label{thm:stratathm}\lremind{thm:stratathm}

The proof is summarized in Section~\ref{pfstratathm}.
In analogy with CW-complexes, we define
an {\em affine cell decomposition} of
a scheme $X$ to be  an affine stratification
$$ X = \bigsqcup_{k} Z_{k} $$
where the $Z_k$'s satisfy (a)--(c) of Theorem~\ref{thm:stratathm}.
The theorem guarantees that such a decomposition exists for any scheme $X$,
with length $\asn X$.

\tpoint{Lemma} {\em Let $\{Y_{k,i}\}$ be an affine stratification of a scheme $X$ and let
 $Z_k := \cup_{i} Y_{k,i}$ be the union of all the affine pieces of index $k$.  
 Then each $Z_k$ is an open dense affine subset of $\overline{Z}_k$,
{\em i.e.}\ $Z_k$ is locally closed and affine. }
 \label{lem:goodunions}\lremind{lem:goodunions}

\bpf
By definition, $Z_k$ is a dense subset of $\overline{Z}_k$.  We will
 see that it is an open subset, and
most importantly, affine. 

Since the affine stratification is finite, we 
have $\overline{\cup_i Y_{k,i}} = \cup_i \overline{Y}_{k,i}$.  
For any distinct $Y_{k,i}$ and $Y_{k,j}$ and any point $y\in \overline{Y}_{k,i}\cap\overline{Y}_{k,j}$, the fact that the $Y$'s are disjoint, 
along with the stratification condition \eqref{bob},
implies that $y$ must be in some $Z_{k'}$ with $k'>k$.  In particular, $y$ is in neither $Y_{k,i}$ nor
$Y_{k,j}$. 

If we let $C_i:=\cup_{j\neq i} \overline{Y}_{k,j}$ be the closed subset consisting of the
closures of other $Y_{k,j}$'s, and $V_i:=X\setminus C_i$ the open complement, then the previous
remark shows that $Y_{k,i}\subseteq V_i$, and therefore that $Z_k\cap V_i= Y_{k,i}$.

Since every locally closed subset is an open subset of its closure,
$Y_{k,i}$ is an open subset of $\overline{Y}_{k,i}\cap
V_i=\overline{Z}_k\cap V_i$.  Since $\overline{Z}_k\cap V_i$ is an
open subset of $\overline{Z}_k$, we see that $Y_{k,i}$ is an open
subset of $\overline{Z}_k$, and therefore that $Z_k = \cup_i Y_{k,i}$
is an open subset of $\overline{Z}_k$.

Let $\widetilde{Z}_k$ be the disjoint union
$$\widetilde{Z}_k = \bigsqcup_{i} Y_{k,i},$$
and $f:\widetilde{Z}_k
\longrightarrow X$ the natural morphism with image $Z_k$.  The map $f$ is
one-to-one on points, and the fact that $Z_k\cap V_i=Y_{i,k}$ for each
$i$ implies that $f$ is a homeomorphism, and in fact an immersion.
Therefore, $\widetilde{Z}_k\cong Z_k$ as schemes, and so $Z_k$ is
affine since $\widetilde{Z}_k$ is.  \epf

\tpoint{Proposition} {\em Let $\{Y_{k,i}\}$ be an affine
  stratification of $X$ of length $m$.  Then there exists an affine
  stratification $\{Y'_{k,j}\}$ of length at most $m$ such that
  the generic points of all components of $X$ are contained in the
  zero stratum $\cup_{j} Y'_{0,j}$ of $\{Y'_{k,j}\}$.}
\label{prop:reorg}\lremind{prop:reorg}

\bpf
We first set $Y'_{0,i}:=Y_{0,i}$ for all valid indices $i$.
Now let $Y_{k,i}$ be any piece of the stratification with $k\geq 1$.  If $Y_{k,i}$ does not
contain the generic point of any component of $X$ then set $Y'_{k,i}:=Y_{k,i}$.
On the other hand, suppose that $Y_{k,i}$ contains $\eta_1,\ldots,\eta_r$ 
where each $\eta_j$ is a generic point of $X$.  
In this case, for each $j \in \{ 1,\ldots,r \}$ 
choose an open affine subset $U_j$ of $Y_{k,i}$ 
containing $\eta_j$ so that $U_j$ intersects no components of $X$ other than
$\clspt{\eta_j}$.   Now set
$Y'_{k,j}:=Y_{k,j}\setminus(\cup_{j=1}^r U_j)$, 
and add the $U_j$ in as elements of the zero stratum, $Y'_{0,i_j}:=U_j$, 
where the $i_j$ are chosen not to conflict with previously existing indices.
It is straightforward to verify that this decomposition satisfies the affine
stratification condition \eqref{bob}.
\epf

\tpoint{Lemma}  {\em
\label{lem:goodclosures}\lremind{lem:goodclosures}Let 
$\{Y_{k,i}\}$ be an affine stratification of $X$ of length $m$.
Then there exists an affine stratification $\{Y'_{k,i}\}$ of length at most
$m$ such
that if $Z'_k:=\cup_{i} Y'_{k,i}$ is the union of all affine pieces of index $k$, then
for any $k$, 
$$\overline{Z}'_k = \cup_{k'\geq k} Z'_{k'}.  $$
}

\bpf
By Proposition~\ref{prop:reorg} we may assume that all the generic points of components
of $X$ occur in the zero stratum of $\{Y_{k,i}\}$, and therefore that 
$\cup_i \overline{Y}_{0,i}=X$.
We now proceed by induction on the length $m$ of the stratification, the case $m=0$ 
being trivial.

Let $U:=\cup_i Y_{0,i}$ be the union of the pieces in the zero stratum, and 
$Z:=\cup_{k\geq 1,i} Y_{k,i}$ the complement. Note that $U$ is open and hence $Z$ is closed
by Lemma~\ref{lem:goodunions}.

The $\{Y_{k,i}\}$ with $k\geq 1$ form 
an affine stratification of $Z$ of length $m-1$
(after reindexing the $k$'s to start with zero). Therefore by induction $Z$ has an
affine stratification of length at most $m-1$ satisfying the hypothesis of the lemma.  
Reindexing the $k$'s again, and adding the $Y_{0,i}$'s as the zero stratum, 
we end up with an affine stratification $\{Y'_{k,i}\}$ of length at most
$m$ 
which also satisfies the hypothesis of the lemma, completing the inductive step. \epf

\tpoint{Corollary} {\em
For any scheme $X$, if $\asn X=m$ then there is
an affine stratification $\{Z_0,\ldots,Z_{m}\}$ of $X$ with 
$\overline{Z}_k=\cup_{k'\geq k} Z_{k'}$ for each $k$, and such that each $Z_k$ is an open
dense affine subset of $\overline{Z}_k$.}\label{cor:goodstrata}\lremind{cor:goodstrata}

\bpf
Combine Lemmas \ref{lem:goodunions} and \ref{lem:goodclosures}. \epf

\tpoint{Corollary} {\em For any scheme $X$, $\asn X \leq\dim X$.}\label{cor:dimboundsasn}
\lremind{cor:dimboundsasn}

\bpf
Let $m=\asn X$ and $\{Z_0,\ldots,Z_{m}\}$ be a stratification as in Corollary
\ref{cor:goodstrata}.
By Corollary~\ref{cor:purecodimone}, each $\overline{Z}_{k+1}$ is of
pure codimension one in $\overline{Z}_k$.  If $\overline{Z}'_{m}$ 
is any irreducible component of $\overline{Z}_{m}$, then that means we can inductively
find a chain of closed irreducible subsets $\overline{Z}'_{m} \subset 
\overline{Z}'_{m-1}\subset\overline{Z}'_{m-2}\subset \cdots \subset 
\overline{Z}'_{1}\subset \overline{Z}'_{0}$, with each $\overline{Z}'_{k}$ an
irreducible component of $\overline{Z}_k$.
Then  $\dim X\geq \dim \overline{Z}'_{m} +m \geq m$.
 \epf
 
 If we assume an additional hypothesis about $X$ or the stratification
 $\{Y_{k,i}\}$, we have slightly stronger results about the
 stratification $\{Z_0,\ldots,Z_m\}$ of
 Corollary~\ref{cor:goodstrata}.

\tpoint{Proposition} {\em If $X$ is an equidimensional scheme, and
  $\{Z_0,\ldots,Z_{m}\}$ the stratification of
  Corollary~\ref{cor:goodstrata}, then we have in addition that
  $\overline{Z}_{k'}$ is of pure codimension $k'-k$ in
  $\overline{Z}_k$ for all $k'\geq k$.  In particular, $\easn X = \asn
  X$.  }\label{equidimX}\lremind{equidimX}

\bpf By Corollary~\ref{cor:purecodimone} each $\overline{Z}_{k+1}$ is
of pure codimension one in $\overline{Z}_k$.  If $\overline{Z}_0=X$ is
equidimensional, then it follows  that each
$\overline{Z}_k$ is equidimensional as well, and from this that
$\overline{Z}_{k'}$ is of pure codimension $k'-k$ in $\overline{Z}_k$
for any $k'\geq k$.  \epf

Even if $X$ is not equidimensional, if the affine stratification $\{Y_{k,i}\}$
satisfies a suitable condition we get a similar good result about the stratification
by the $Z_k$'s.

\tpoint{Proposition} {\em Let $\{Y_{k,i}\}$ be an affine stratification of a
  scheme $X$ with each $Y_{k,i}$ of pure codimension $k$ in $X$.  Let
  $Z_k := \cup_{i} Y_{k,i}$ be the union of all the affine pieces of
  codimension $k$.
  Then\label{stratagivecells}\lremind{stratagivecells}
\begin{itemize}
\item[(i)] For $k'\geq k$, $\overline{Z}_{k'}$ is of pure codimension 
$k'-k$ in $\overline{Z}_{k}$; in particular, $\overline{Z}_{k'}\subseteq \overline{Z}_k$.
\item[(ii)]  $\overline{Z}_k = \cup_{k'\geq k} Z_{k'}$.
\end{itemize}
}

\bpf Since the decomposition is finite, the irreducible
components of $\overline{Z}_{m}$ are all of the form
$\overline{W}_{m}$ with $W_{m}$ an irreducible component of some
$Y_{m,j}$.

We prove (i) by induction on $k$.  For $k=0$ the result is obvious,
since $\overline{Z}_0=X$, and $\overline{Z}_{k'}$ is of pure
codimension $k'$ in $X$. So assume that $k>0$ and that (i) is true for
$k-1$.

Let $\overline{W}_{k'}$ by any irreducible component of
$\overline{Z}_{k'}$ with $k'\geq k$.  By the induction hypothesis,
$\overline{W}_{k'}\subset \overline{Z}_{k-1}$, and is of codimension
$k'-k+1$ in $\overline{Z}_{k-1}$.  Let $T_{k-1}$ be any irreducible
component of $Z_{k-1}$ whose closure contains $\overline{W}_{k'}$ and
such that $\codim(\overline{W}_{k'},\overline{T}_{k-1})= k'-k+1$.
Lemma~\ref{lem:goodunions} gives us that $Z_{k-1}$ is affine, and
therefore $T_{k-1}$ is affine also.

By  Lemma~\ref{codimone}, the closed set $\overline{T}_{k-1}\setminus T_{k-1}$ 
has codimension one in $\overline{T}_{k-1}$.
Let $\eta_k$ be the generic point of any component of $\overline{T}_{k-1}\setminus T_{k-1}$ containing $W_{k'}$; one exists by our
choice of $T_{k-1}$.  Since 
$\codim(\overline{W}_{k'},\overline{T}_{k-1})=k'-k+1$
and
$\codim(\clspt{\eta_k},\overline{T}_{k-1})=1$,
we have $\codim(W_{k'},\clspt{\eta_k})=k'-k$.

The $Z_{m}$'s partition $X$, and so $\eta_k$ must be in exactly one
$Z_{m}$.  We cannot have $m\leq k-1$, since that would
contradict the stratification condition.  We cannot have $m\geq
k+1$, since this would contradict
$\codim(\overline{Z}_{m},\overline{Z}_{k-1})=m-k+1$, which
holds by the induction hypothesis.  Therefore $\eta_k$ is in $Z_k$,
and so $\overline{W}_{k-1}\subset \overline{Z}_k$.

We have already seen that $\codim(W_{k'},\clspt{\eta_k})=k'-k$.  Since 
$$\codim \left( \overline{W}_{k'},\overline{Z}_k \right) =
\sup_{i} \left( \codim \left( \overline{W}_{k'},\overline{Y}_{k,i} \right)
\right),$$
with
$\overline{Y}_{k,i}$ running over the components of $\overline{Z}_k$,
we get $\codim(\overline{W}_{k'},\overline{Z}_k)\geq k'-k.$ But for
any three closed schemes $\overline{W}$, $\overline{Z}$, and $X$ with
$\overline{W}\subseteq\overline{Z}\subseteq X$, we always have 
$$\codim \left( \overline{W},\overline{Z} \right)+\codim \left( \overline{Z},{X} \right)\leq \codim \left( \overline{W},X \right).$$
Since the codimensions of $\overline{W}_{k'}$ in $X$, and $\overline{Z}_k$ in $X$ are $k'$ and $k$ by hypothesis, this gives
$\codim(\overline{W}_{k'},\overline{Z}_k)\leq k'-k$, and hence
$\codim(\overline{W}_{k'},\overline{Z}_k)= k'-k.$ 
Therefore $\overline{Z}_{k'}$ is contained in $\overline{Z}_k$, and is of pure codimension $k'-k$, completing the inductive step for (i).

To prove (ii), the stratification condition gives 
$\overline{Z}_k\subseteq \cup_{k'\geq k} Z_{k'}$, while part (i) above 
gives the opposite inclusion. \epf

\epoint{Proof of Theorem~\ref{thm:stratathm}}
\label{pfstratathm}\lremind{pfstratathm}Part 
(a) is Lemma~\ref{lem:goodclosures}, part (b) Lemma~\ref{lem:goodunions}, part (c)
Corollary~\ref{cor:purecodimone}, part (d) Proposition~\ref{equidimX}, and part (e) Proposition
\ref{stratagivecells} and Lemma~\ref{lem:goodunions} again. \epf

\section{Topological consequences of bounded affine stratification number}
\label{topcon}

We now describe the topological consequences of bounded $\asn$, in
particular: relation to dimension (Section~\ref{reldim}), affine
covering number (Section~\ref{relacn}), cohomological dimension (for
quasicoherent sheaves, Section~\ref{relcd}, as well as constructible and
$\ell$-adic sheaves, Section~\ref{cv}), dimension of largest proper
subscheme (Section~\ref{relpsv}), and homotopy type (Section~\ref{relht}).

\bpoint{Relation to dimension}\label{reldim}

\tpoint{Proposition}  {\em $\asn X \leq \dim X$.  If one
top dimensional component of $X$ is proper, then equality holds.}
\label{asnleqdim}\lremind{asnleqdim}

The first statement is Corollary~\ref{cor:dimboundsasn}.
The second statement follows from
Proposition~\ref{cdleqasn} ($\cd \leq \asn$) and the following theorem,
first conjectured by Lichtenbaum.

\tpoint{Theorem (Grothendieck \cite[6.9]{lc}, Kleiman \cite[Main Theorem]{kleiman})}
{\em If $d=\dim X$, then $\cd X=d$ if and only if at least one
  $d$-dimensional component of $X$ is proper.}\label{lichtenbaum}
\lremind{lichtenbaum} 

\epoint{Example: All values between $0$ and $\dim X$ are possible} Let
$X_k = \PP^n\setminus\{\mbox{$(n-k-1)$-plane)}\}$, for $k$ between
$0$ and $n-1$.  Then clearly $\cd X_k=k$ and $\asn X_k  \leq k$.  
We will see
that $\cd \leq \asn$ (Proposition~\ref{cdleqasn}), from which the result follows.

\bpoint{Relation to affine covering number}
\label{relacn}\lremind{relacn}Recall that the affine covering number
$\acn X$ of a scheme $X$ is the minimal number of affine open subsets
required to cover $X$, minus $1$.  The invariant $\acn$ does not
obviously behave as well as $\asn$ with respect to products ({\em
  cf.}\ Lemma~\ref{triviallemma}(c)); it also is not bounded
by dimension (Example~\ref{jason} below).

The argument of
\eqref{exists} gives the following.

\tpoint{Proposition} {\em  $\asn X \leq \acn X$.} \label{asnleqacn}
\lremind{asnleqacn}

\epoint{Example} 
\label{k3}\lremind{k3}In general, $\acn X \neq \asn X$.  As an
example, let $X$ be a complex 
K3 surface with Picard rank $1$, minus a very
general point.  Then $\acn X= 2$: if for every point $p$ of $X$, $\acn
(X-p) =1$, then (given the hypothesis that the Picard rank is $1$) any
two points of $X$ are equivalent in $A_0(Y)$ (with $\Q$-coefficients),
contradicting Mumford's theorem that $A_0(Y)$ is not countably
generated \cite{theword}.  Example~\ref{jason} below gives another example
(in light of Proposition~\ref{asnleqdim}).

\epoint{Example: $\acn X$ may be larger than $\dim X$}
\label{jason}\lremind{Jason}When $X$ is quasiprojective, $\acn X \leq
\dim X$.  ({\em Reason:} Let $\overline{X}$ be a projective
compactification such that the complement $\overline{X}\setminus X$ is
a {Cartier} divisor $D$.  Consider an embedding $\overline{X}
\hookrightarrow \proj^n$ and let $H_0,\ldots,H_{\dim X}$ be
hypersurfaces so that $\overline{X} \cap H_0 \cap \cdots \cap H_{\dim
  X} = \emptyset$.  Then the $\overline{U_i}:=
\overline{X}\setminus(\overline{X}\cap H_i)$ form an affine cover of
$\overline{X}$.  We conclude using Lemma~\ref{minuscartier}.)  

However, the following example, due to J. Starr, shows that $\acn X$
may be greater than $\dim X$.  Given any $n$, we describe a reducible,
reduced threefold that requires at least $n$ affine open sets to cover
it.  Recall Hironaka's example ({\em e.g.}\ \cite[Example~B.3.4.1]{h}) of
a nonsingular proper nonprojective threefold $X$.
Nonprojectivity is shown by exhibiting two curves $\ell$ and $m$
whose sum is numerically trivial.  Hence no affine open set can meet
both $\ell$ and $m$; otherwise its complement would be a divisor
(Lemma~\ref{codimone}), hence Cartier (as $X$ is nonsingular), which
meets both $\ell$ and $m$ positively.  Now choose points $p$ and $\ell$
and $q$ on $m$.  Consider $\binom n 2$ copies of $(X, p, q)$, corresponding
to ordered pairs $(i,j)$ ($1 \leq i<j \leq n$); call these
copies $(X_{ij}, p_{ij}, q_{ij})$.  Let $r_1$, \dots, $r_n$ be copies
of a reduced point.  Glue $r_i$ to $p_{ji}$ and $q_{ik}$.
Then no affine open can contain both $r_i$ and $r_j$ for $i < j$ (by
considering $X_{ij}$).

\epoint{Example: a family of integral threefolds with arbitrary high
  affine covering number (and no smooth embeddings)}
\label{nosmooth}This leads to an example of an integral (but singular)
threefold that requires at least $n$ affine open sets to cover it.
(Question: Is there a family of nonsingular irreducible varieties of
fixed dimension with unbounded affine covering number?)  Our example
will be a blow-up of $\proj^3$.  Choose $n$ curves $C_1$, \dots, $C_n$
in $\proj^2 \subset \proj^3$ that meet in $n$ simple $n$-fold points
$p_1$, \dots, $p_n$ (and possibly elsewhere).  Away from $p_1$, \dots,
$p_n$ blow up $C_1$, \dots, $C_n$ in some arbitrary order.  In a
neighborhood of $p_i$ (not containing any other intersection of the
$C_j$) blow up $C_i$ first (giving a smooth threefold) and then blow
up the local complete intersection $\cup_{j \neq i} C_j$ (or more
precisely, the proper transform thereof), giving a threefold with a
single singularity (call it $q_i$).  The preimage of $p_i$ is the
union of two $\proj^1$'s, one arising from the exceptional divisor of
$C_i$ (call it $\ell_i$), and one from the exceptional divisor of
$\cup_{j \neq i} C_j$; they meet at $q_i$.  By Hironaka's argument,
$\ell_i + \ell_j$ is numerically trivial for all $i \neq j$.  Then no
affine open $U$ can contain both $q_i$ and $q_j$: the complement of
$U$ would be a divisor, meeting $\ell_i$ and $\ell_j$ properly and at smooth
points of our threefold (i.e.\ not at $q_i$ and $q_j$), and the same
contradiction applies.

This is also an example of a scheme which cannot be embedded in any
smooth scheme, or indeed algebraic space.  (Earlier examples are the
topic of papers of Horrocks \cite{ho} and Nori \cite{n}.)  If $X
\hookrightarrow W$ with $W$ smooth, then for any divisor $D$
(automatically Cartier) on $W$, $D\cdot \ell_i + D\cdot \ell_j=0$ for
all $i,j$.  If $n\geq 3$, this implies that $D\cdot\ell_i=0$ for all
$i$.  But for any affine open set $U$ of $W$ with $U\cap
\ell_i\neq\emptyset$, the complement $D=W\setminus U$ would intersect
$\ell_i$ properly, giving us the contradiction $D\cdot \ell_i>0$.
Hence no such embedding is possible.

By combining Proposition~\ref{asnleqdim} with
Proposition~\ref{asnleqacn}, we obtain the following.

\tpoint{Proposition}  {\em If $X$ is proper, then $\acn X \geq \dim
  X$.}

We note that this also follows from Theorem
\ref{lichtenbaum}.
Example~\ref{k3} shows that it is {\em not} true that $\acn X = \dim
X$ if and only if $X$ is proper, even for quasiprojective $X$.

\bpoint{Relation to cohomological dimension}\label{relcd}
Just as the dimensions of the cells in a CW-complex bounds the topological (co)homology, the
length of the stratification into affine cells bounds the quasicoherent sheaf cohomology.

\tpoint{Proposition} {\em $\cd X\leq \asn X$.} \label{cdleqasn}\lremind{cdleqasn}

\bpf We prove the result by induction on $\asn X$.  It is clear for
$\asn X=0$, so assume that $m:=\asn X>0$ and that the result is proven
for all schemes $Z$ with $\asn Z <\asn X$.  Let $\{Z_0,\ldots,Z_{m}\}$
be an affine cell decomposition given by Theorem~\ref{thm:stratathm}.
Set $Z:=X\setminus Z_0 = \overline{Z}_1 = \cup_{k\geq 1} Z_k$.

By Theorem~\ref{thm:stratathm}(b) $Z_0$ is an open dense
affine subset of $X$, so by Corollary~\ref{inductiveboundingstep},
$\cd X\leq \cd Z+1$.
Next, $Z = \sqcup_{k=1}^{m-1} Z_k$ is (after reindexing) an affine stratification of $Z$ of length
$m-1$, so $\asn Z \leq \asn X-1$.
Finally, by the inductive hypothesis,
$\cd Z\leq\asn Z$. 
Combining these three inequalities gives $\cd X\leq\asn X$,
completing the inductive step. \epf

We remark in passing that by combining Proposition~\ref{cdleqasn}
with Corollary~\ref{cor:dimboundsasn} we obtain another proof of
Grothendieck's dimensional vanishing theorem (\cite[Theorem\ 
3.6.5]{tohoku}, \cite[Theorem\ III.2.7]{h}).

We conclude with an obvious result.

\tpoint{Proposition}  {\em $\cd X = 0$ if and only if $\asn X = 0$ if and
only if $\acn X = 0$.}

\bpf Each of the three is true if and only if $X$ is affine (the first by
Serre's criterion for affineness). \epf

\bpoint{Relation to dimension of largest complete subscheme}
\label{relpsv}\lremind{relpsv}Motivated 
by Diaz' theorem \cite{diaz}, let $\psv X $ be
the largest dimension of a proper closed subscheme of $X$.
If $Z$ is a proper closed subscheme of $X$ (with inclusion
$j:Z\hookrightarrow X$), and if $\cF$ is a quasicoherent sheaf on $Z$,
then $j_{*}\cF$ is a quasicoherent sheaf on $X$, and
$H^{i}(X,j_{*}\cF) = H^i(Z,\cF)$ for all $i$.  
By Theorem~\ref{lichtenbaum} we can find a
quasicoherent sheaf $\Fsh$ on $Z$ with $H^{\dim Z}(Z,\Fsh)\neq 0$, and so
this gives $\psv X \leq \cd X.$  Hence by Proposition~\ref{cdleqasn},
$$\psv X \leq \asn X.$$

\bpoint{Relation to cohomological vanishing for constructible and
  $\ell$-adic sheaves} \label{cv}\lremind{cv}In this section all
notions related to sheaves (including stalks, pushforwards, and
cohomology groups) are with respect to the \'{e}tale topology.  For
instance, ``sheaf on $X$'' means ``sheaf on $X$ in the \'etale
topology''.

To show how $\asn$ implies cohomological vanishing for constructible and
$\ell$-adic sheaves (Corollary~\ref{fly}), we first
recall a theorem and some notation of Artin.
For any (\'{e}tale) sheaf $\tsF$ of abelian groups on $X$, let
$$\dimsup(\tsF): =\sup\left\{{\dim \left( \clspt{x} \right)\,\,|\,\, x\in X,\,\,\tsF_{\overline{x}}\neq 0}\right\} $$ 
be the dimension of the support of $\tsF$.

\tpoint{Artin's Theorem \cite[Theorem~3.1]{artin}} {\em
Let $f:X\longrightarrow Y$ be an affine morphism of schemes of finite type
over a field $k$, and $\tsF$ a torsion sheaf ({\em i.e.}\ sheaf of torsion
groups) on $X$.   Then 
$\dimsup(R^qf_{*}\tsF)\leq \dimsup(\tsF)-q$ for all $q\geq 0$.
}\label{thm:artin}\lremind{thm:artin}

We will apply Artin's Theorem in the following form:

\tpoint{Proposition} {\em Suppose that $X$ is a scheme, $U$ an affine
  open subset of $X$, and $Z:=X\setminus U$ the complement.  Then for
  any torsion sheaf $\tsF$ on $X$,
$$\dimsup(\ch_{\et,Z}^q(\tsF))\leq \dimsup(\tsF) - q+1.$$}\label{prop:localtopbound}\lremind{prop:localtopbound}

Here the $\ch_{\et,Z}^q(\tsF)$ are the local cohomology sheaves in the
\'{e}tale topology.  The usual excision and spectral sequences for
local cohomology remain true in the \'etale setting, see \cite[Sec.
6]{theverd}.

\bpf
If  $i:U\hookrightarrow X$ is the inclusion, then for any sheaf $\tsF$
of abelian groups on $X$ we have the exact sequence \cite[Proposition~6.5]{theverd}
$$ 0 \longrightarrow \ch_{\et,Z}^0(\tsF) \longrightarrow \tsF \longrightarrow i_{*}(\tsF|_U)
\longrightarrow \ch_{\et,Z}^1(\tsF)   \longrightarrow 0,$$
as well as isomorphisms
$$ \ch_{\et,Z}^{q}(\tsF) \cong R^{q-1}i_{*}(\tsF|_U) \,\, \mbox{for all $q\geq 2$}.$$

If $q\geq 2$ the proposition then follows from the above isomorphism
and Artin's Theorem~\ref{thm:artin} applied to the inclusion morphism
$i$, which is affine since $U$ is.

If $q=1$
it is enough to bound $\dimsup(i_{*}(\tsF|_U))$, since $\ch_{\et,Z}^1(\tsF)$ is a quotient of $i_{*}(\tsF|_U)$.
The points $x\in X$ where $(i_{*}(\tsF|_U))_{\overline{x}}\neq 0$ are the points $x\in U$ with $\tsF_{\overline{x}}\neq 0$ and
points $x\in Z$ such that there exists a point $x'\in U$, $x\in\clspt{x'}$ with $\tsF_{\overline{x'}}\neq 0$.  In particular,
the support of $i_{*}(\tsF|_U)$ is contained in the support of $\tsF$, so $\dimsup(\ch_{\et,Z}^1(\tsF))\leq\dimsup(\tsF)$,
which is exactly the statement of the proposition when $q=1$.

If $q=0$ note that  $\dimsup(\ch_{\et,Z}^0(\tsF))\leq \dimsup(\tsF)$ since
$\ch_{\et,Z}^0(\tsF)$ is a subsheaf of $\tsF$, while the proposition only claims the weaker bound
  $\dimsup(\ch_{\et,Z}^0(\tsF))\leq \dimsup(\tsF) +1$.
\epf
 
\tpoint{Lemma} {\em If $\tsF$ is a torsion sheaf on $X$, then $H_\et^n(X,\tsF)=0$ for all
$n>\dimsup(\tsF)+\asn X.$ }\label{lem:torsionbound}\lremind{lem:torsionbound}

\bpf We show the result by induction on $\asn X$, the case $\asn X=0$
being Artin's Theorem~\ref{thm:artin} again.  Let
$\{Z_0,\ldots,Z_{\asn X}\}$ be an affine cell decomposition of $X$ (as
given by Theorem \ref{thm:stratathm}).  Set
$Z:=X\setminus Z_0=\cup_{k\geq 1}Z_k$.

We have $H_{\et}^n(Z_0,\tsF|_{Z_0})=0$ for all $n>\dimsup(\tsF|_{Z_0})$ by
Artin's Theorem, and since $\dimsup(\tsF|_{Z_0})\leq\dimsup(\tsF)$ the
excision sequence (equation~\eqref{excision} holds in this context,
\cite[(6.5.3)]{theverd}) shows that $H_{\et}^n(X,\tsF)$ is a quotient
of $H_{\et,Z}^n(\tsF)$ for all $n>\dimsup(\tsF)$. It is therefore
enough to show that $H_{\et,Z}^n(\tsF)=0$ for $n>\dimsup(\tsF)+\asn
X$.

We can compute $H^n_{\et,Z}(\tsF)$ by a spectral sequence with $E_2^{pq}$ term
$H^p_{\et}(X,\ch_{\et,Z}^q(\tsF))$ (\cite[Proposition 6.4]{theverd}). We have $H^p_{\et}(X,\ch_{\et,Z}^q(\tsF))=H^p_{\et}(Z,\ch_{\et,Z}^q(\tsF))$ since 
$\ch_{\et,Z}^{q}(\tsF)$ is supported on $Z$.
By Proposition~\ref{prop:localtopbound} we have $\dimsup(\ch_{\et,Z}^q(\tsF))\leq \dimsup(\tsF)-q+1$.
Since $\asn Z <\asn X$ we can apply the inductive hypothesis to conclude that
$H^p_{\et}(Z,\ch_{\et,Z}^q(\tsF))=0$ for $p>\dimsup(\tsF)-q+1+\asn Z $, or
$p+q>\dimsup(\tsF)+\asn Z +1$.
Again using $\asn Z <\asn X$, this gives $H^n_{\et,Z}(\tsF)=0$ for $n>\dimsup(\tsF)+\asn X$.
\epf

\tpoint{Corollary}{\em 
\begin{itemize}
\item[(a)] If $\tsF$ is a torsion sheaf, then $H^n_{\et}(X,\tsF)=0$ for all $n>\dim X +\asn X$.
\item[(b)] If $\tsF$ is a constructible sheaf, then $H^n_{\et}(X,\tsF)=0$ for all $n>\dim X +\asn X$.
\item[(c)] If $\cF_{\ell}$ is an $\ell$-adic sheaf on $X$, then $H^n_{\et}(X,\cF_{\ell})=0$ for all 
$n>\dim X +\asn X$.
\end{itemize}}
\label{fly}

\bpf
{\em (a)} \lremind{fly}Clearly we have $\dimsup(\tsF)\leq\dim X$. 
{\em (b)} A constructible sheaf is a special case of a torsion sheaf (compare \cite[Proposition 1.2(ii)]{artin2} with
\cite[Def. 2.3]{artin2}).
{\em (c)} follows from (a).
\epf

\bpoint{Relation to homotopy type}
\label{relht}\lremind{relht}We 
expect that the affine stratification number bounds the homotopy type
as follows.

\tpoint{Conjecture} {\em If the base field is $\C$, then $X$ has the
  homotopy type of a finite complex of dimension at most $\asn X+\dim
  X$.}  \label{homtype} \lremind{homtype}

\section{Applications to moduli spaces of curves}
\label{ms}\lremind{ms}One  
motivation for the definition of affine stratification
number is the study of the moduli space of curves, and certain
geometrically important open subsets.    We will use
Definition~\ref{defeasn} (which we may, by Proposition~\ref{equidimX}).

\epoint{Preliminary aside: the affine stratification number of
  Deligne-Mumford stacks} As we have only defined the affine stratification
number of schemes, throughout this section, we will work with coarse
moduli space of curves.  One should presumably work instead with a
more general definition for Deligne-Mumford stacks.  One possible
definition is to replace the notion of ``affine'' in the definition of
affine stratification number with that of a Deligne-Mumford stack that
has a surjective finite flat cover by an affine scheme (see Proposition\ 
\ref{macbeth}).

\point
Recall the following question of Looijenga's.

\tpoint{Conjecture (Looijenga)}
{\em (a)  $\acn M_g \leq g-2$ for $g \geq 2$.
(b) More generally, $\acn M_{g,n} \leq g-1 - \de_{n,0}$ whenever
$g>0$, $(g,n) \neq (1,0)$.}

The case $n=1$ of (b) implies the cases $n>1$, as the morphism $M_{g,n+1}
\rightarrow M_{g,n}$ is affine for $n \geq 1$.

This suggests the following, weaker conjecture, which is
straightforward to verify for small $(g,n)$ (using Proposition~\ref{macbeth} 
judiciously).  We are currently
pursuing a program to prove this (work in progress).

\tpoint{Conjecture (Looijenga \cite[p.~112, Problem~6.5]{hl})} 
{\em $\asn  M_g  \leq g-2$ for $g \geq 2$.}
\label{key} \lremind{key}

From this statement (and properties of $\asn$), we obtain a number
of consequences.  

\tpoint{Proposition (Looijenga \cite[p.~112]{hl})} 
{\em Conjecture~\ref{key} implies that
$\asn M_{g,n} \leq g-1 -\de_{n,0}$ 
 whenever
$g>0$, $(g,n) \neq (1,0)$.} \label{key1}
\lremind{key1}

\bpf As $M_{g,n+1} \rightarrow M_{g,n}$ is affine for $n \geq 1$, it
suffices to prove the result for $M_{0,3}$ and $M_{g,1}$ with $g>0$. 
The cases $g=0$ and $g=1$ are immediate.  For $g>1$, let
$D$ be a multisection of $M_{g,1} \rightarrow M_g$ ({\em e.g.}\  a suitable
Weierstrass divisor).  Then the morphisms $D \rightarrow M_g$ and $(M_{g,1}\setminus D)
\rightarrow M_g$ are affine and surjective, so pulling back the
affine stratification of $M_g$ to $M_{g,1}$ and intersecting with
$(M_{g,1} \setminus D) \sqcup D$ yields the desired affine stratification of
$M_{g,1}$. \epf

Examination of small genus cases suggests 
the following refinement of Conjecture~\ref{key}.

\tpoint{Conjecture} {\em There is an affine stratification of $M_{g',
    n'}$ preserved by the symmetric group acting on the $n'$ points.
  The induced decomposition of $\Mbar_{g,n}$ is a stratification.}
\label{key2} \lremind{key2}

This leads to a bound on the affine stratification
number of the open subset $\Mbar_{g,n}^{\leq k}$, corresponding to
stable $n$-pointed genus $g$ curves with at most $k$ genus $0$
components, defined in \cite[Section~4]{gv}.

\tpoint{Proposition} {\em Conjecture~\ref{key2} implies
that $\asn  \Mbar_{g,n}^{\leq k} \leq g-1+k$ for all $g>0$, 
$n>0$.} \label{response} \lremind{response}

This is more evidence of the relevance of this strange filtration of
the moduli space of curves.  In particular, compare this to Theorem~$\star$ of \cite{gv}, that the tautological ring of $\Mbar_{g,n}^{\leq
k}$ vanishes in codimension greater than $g-1+k$.  
(In \cite[p.~3]{gv}, Looijenga asks precisely this question,
with $\asn$ replaced by $\acn$.)

\bpf
We show that dimension of any stratum
of $\Mbar_{g,n}$ appearing in $\Mbar_{g,n}^{\leq k}$ is 
at least 
$$3g-3+n-(g-1+k) = 2g-2+n-k.$$
This is true for  strata  in $M_{g,n}$ by Proposition~\ref{key1}.
Consider any other boundary
stratum, say with $j$ rational components ($j \leq k$) with $m_1,
\dots, m_j$ special points respectively; and $s$ other components,
with genus $g_1$, \dots, $g_s$ and $n_1$, \dots, $n_s$ special points
respectively.  By Proposition~\ref{macbeth}, it
suffices to pass to the finite \'etale cover that is isomorphic to
$$\prod^j_{i=1} M_{0,m_i} \times \prod_{i=1}^s M_{g_i, n_i}.$$
By Proposition~\ref{key1} (and using
$M_{0,m_i}$ affine), we can decompose this space into affine sets of dimension
at least
$$
 \sum_{i=1}^j (m_i-3)  + \sum_{i=1}^s ( 3g_i-3+n_i - (g_i-1)).
$$
Now $\sum_{i=1}^j (m_i-2) + \sum_{i=1}^s (2g_i-2+n_i) = 2g-2+n$, so
each affine set has dimension at least
$$
2g-2+n - j
$$
and thus codimension in $\Mbar_{g,n}^{\geq k}$ at most
$g-1+j \leq g-1+k$ as desired.
\epf

This leads to bounds on other spaces of interest.  Let $M_{g,n}^{ct}$
be the open subset of $\Mbar_{g,n}$ corresponding to curves of compact
type ({\em i.e.}\ with compact Jacobian, or equivalently with dual
graph containing no loops).  Let $M_{g,n}^{rt}$ be the open subset
corresponding to curves with rational tails ({\em i.e.}\ with a
component a smooth genus $g$ curve, or equivalently with dual graph
with a genus $g$ vertex).

\tpoint{Corollary}
{\em 
Conjecture~\ref{key2} implies
that $\asn  M_{g,n}^{ct}  \leq 2g-3+n$ and 
$\asn  M_{g,n}^{rt}  \leq g+n-2$ for $g>0$, $n>0$.
}

\bpf $M_{g,n}^{ct}$ is obtained by removing boundary strata from
$\Mbar^{\leq g+n-2}_{g,n}$.  $M_{g,n}^{rt}$ is obtained by removing
boundary strata from $\Mbar^{\leq n-1}_{g,n}$. 
\epf

\tpoint{Corollary}  {\em Conjectures~\ref{homtype} and \ref{key2} imply
that $\Mbar_{g,n}^{\leq k}$ (resp.\ $M_{g,n}^{ct}$, $M_{g,n}^{rt}$)
has the homotopy type of a finite complex of dimension at most
$4g-4+n+k$ (resp.\ $5g-6+2n$, $4g-5+2n$).}
\label{homMgn} \lremind{homMgn}

} % end of parskip; it started just before the introduction

\end{document}